\documentclass[11pt,english]{smfart}
\usepackage[applemac]{inputenc}
\def\today{mars 2008} 
\usepackage{amsmath,amsfonts,amsthm,amssymb,amscd}
\usepackage[francais,english]{babel}
\renewcommand{\Im}{\mathop{\rm Im}\nolimits}

\newcommand{\be}{\begin{equation}}
\newcommand{\ee}{\end{equation}}
\newcommand{\bi}{\begin{itemize}}
\newcommand{\ei}{\end{itemize}}

\theoremstyle{plain} \newtheorem{theorem}{Theorem}[section]
\newtheorem{lemma}[theorem]{Lemma}
\newtheorem{proposition}[theorem]{Proposition}

 \theoremstyle{remark}

\newcommand{\R}{{\mathbb R}}

\newcommand{\Z}{{\mathbb Z}}
\newcommand{\N}{{\mathbb N}}
\newcommand{\e}{{\varepsilon}}

\newcommand{\Ac}{{\mathcal A}} 
\newcommand{\Mc}{{\mathcal M}}

\newcommand{\Om}{\Omega}

\def\sleq{\leq\kern-6pt \cdot\null\hskip4pt}

\numberwithin{equation}{section}

\begin{document}

\author{ Beno\^it Gr\'ebert, Carlos Villegas-Blas} \title[On energy exchanges in NLS]{On the energy exchange between resonant modes in  nonlinear Schr\"odinger equations}
\alttitle{Echange d'énergie entre modes résonants dans une équation de Schr\"odinger non linéaire cubique.}
\date{\today}
\thanks{\noindent B. Grébert  supported in part by the grant ANR-06-BLAN-0063.\\ C. Villegas-Blas supported in part by PAPIIT-UNAM IN106106-2}
\begin{abstract}
We consider the nonlinear Schr\"odinger equation
$$ i\psi_t= -\psi_{xx}\pm 2\cos 2x \ |\psi|^2\psi,\quad x\in S^1,\ t\in \R$$
and we prove that the solution of this equation, with small initial datum $\psi_0=\e (\cos x+\sin x)$, will periodically exchange energy between the Fourier modes $e^{ix}$ and $e^{-ix}$. This beating effect is described up to time of order $\e^{-9/4}$ while the frequency is of order $\e^2$. We also discuss some generalizations. 

    \end{abstract}

\begin{altabstract} 
 Nous considérons l'équation de Sch\"odinger non linéaire 
 $$ i\psi_t= -\psi_{xx}\pm 2\cos 2x \ |\psi|^2\psi,\quad x\in S^1,\ t\in \R$$
 et nous montrons la solution de cette équation ayant pour donnée initiale $\psi_0=\e (\cos x+\sin x)$ avec $\e$ petit, va échanger périodiquement de l'énergie entre les modes de Fourier $e^{ix}$ et $e^{-ix}$. Cet effet de battement, dont la période est de l'ordre de $\e^{-2}$, est mis en évidence pour des temps de l'ordre de $\e^{-9/4}$. Nous présentons aussi quelques généralisations. 
\end{altabstract}    
\altkeywords{Forme Normale, Equation de Schr\"odinger non linéaire, 
résonances, échange d'énergie}\frontmatter
\maketitle

\setcounter{section}{0}

\section{Introduction}
Let us consider the non-linear Schr\"{o}dinger equation (NLS) on the circle $S^1$
\begin{equation}\label{NLSE}
\imath \psi_t = - \psi_{xx}+ V*\psi +g(x,\psi,\bar{\psi})
\end{equation}
where $V*\Psi$ denotes the convolution function between a potential $V:S^1\mapsto\mathbb{R}$ and the function $\Psi$. The nonlinear term $g$ is  Hamiltonian in the sense that $g= \partial_{\bar \psi}G$ with $G$ analytic with respect to its three variables and $G(x,z,\bar z)\in \R$. We further assume that $G$ is globally at least of order three in $(z,\bar z)$ at the origin in such a way  $g$ is effectively a nonlinear term.\\
 The Fourier basis
$\exp(i{j}x)$, $j\in\mathbb{Z}$ provides an orthonormal basis for $L^2(S^1)$ in which the linear operator $A=- \frac{\partial^2}{\partial{x}^2} + V*$ is diagonal.  The corresponding eigenvalues of $A$ are
given by the real numbers
 $\omega_j= j^2 + \hat{V}(j)$,  $j\in\mathbb{Z}$, where $\hat{V}(j)$ are the Fourier coefficients of $V$:
\begin{equation}\label{Foucoeff}
\hat{V}(j) = \int_{-\pi}^{\pi}V(x)\exp(-\imath{j}x)dx
\end{equation}
where $dx$ denotes the normalized Lebesgue measure on $S^1$.

Under a non-resonant condition on the frequencies $\omega_j$,
it has been established by D. Bambusi, B. Grebert  \cite{BG06} (see also \cite{Gre07})
that, given a number $r\geq 3$, for small initial data in $H^s$ norm\footnote{Here $H^s$ denotes the standard Sobolev Hilbert space on $S^1$ and $||\cdot ||_s$ its associated norm.}, say $||\psi_0||_s=\e$ with $s$ large enough, the solution to the NLS equation (\ref{NLSE}) remains small  in the same Sobolev norm, $||\psi(t)||_s\leq 2\e$ for a large period of time, $|t|\leq \e^{-r}$. Furthermore  the actions
$I_j\equiv\xi_j\eta_j$, $j\in\mathbb{Z}$, are almost invariant during the same period of time. In \cite{BG06}, it is also proved that in the more natural case of multiplicative potential,  
$$\imath \psi_t = - \psi_{xx}+ V\psi +g(x,\psi,\bar{\psi}),$$
which corresponds to an asymptotically resonant case, 
$\omega_j\sim\omega_{-j}$ when $|j|\to \infty$,  we can generically impose non resonant conditions on  the frequencies $\lbrace\omega_{|j|}\rbrace$ in such a way that the generalized actions $J_j=I_j+I_{-j}$ are almost conserved quantities. A priori nothing prevent  $I_j$ and
$I_{-j}$ from interacting. In this article,  we exhibit nonlinearity $g$ that incites this interaction and especially a beating effect.  This is a nonlinear effect which is a consequence of the  resonances of the linear part. We will see that not all nonlinearities incite the beating. For instance if $g$ does not depend on $x$ there is no beating for a long time. Typically, to obtain an interaction between the mode $j$ and the mode $-j$, the nonlinearity $g$ must contain oscillations of frequency $2j$ or a multiple of $2j$, i.e. $g$ must depend on $x$ as $\cos 2kjx$ or $\sin 2kjx$ for some $k\geq 1$. In what follows, we will focus on the totally resonant case $V=0$ and a cubic nonlinearity $g=2\cos 2x \ |\psi|^2\psi$. Precisely
we consider the following Cauchy problem
\be \label{cauchy}
\left\{\begin{array}{ll}i\psi_t & = -\psi_{xx}\pm 2\cos 2x \ |\psi|^2\psi,\quad x\in S^1,\ t\in \R\\
\psi(0,x) &=\e (\cos x +\sin x)
\end{array}\right.
\ee
where $\e$ is a small parameter. By classical arguments based on the conservation of the energy and Sobolev embeddings, this problem has a unique global  solution in all the Sobolev spaces $H^s$ for $s\geq 1$ and even for $s\geq 0$ using more refined technics (see for instance \cite{Bou99} chap. 5).  Our result precises the behavior of the solution for time of order $\e^{-9/4}$. We write the solution $\psi$ in Fourier modes: $\psi(t,x)=\sum_{k\in \Z}\hat \psi_k(t) e^{ikx}$.
\begin{theorem}\label{main}
For $\e$ small enough we have for $|t|\leq \e^{-9/4}$
\be \left\{\begin{array}{ll}
|\hat\psi_1(t)|^2+|\hat\psi_{-1}(t)|^2&=\e^2+O(\e^{3}), \\
|\hat\psi_1(t)|^2-|\hat\psi_{-1}(t)|^2&=\pm\e^2\sin 2\e^2t+O(\e^{9/4}).\end{array}\right.
\ee
\end{theorem}
In other words, the first estimate says that all the energy remains concentrated on the two Fourier modes $+1$ and $-1$ and the second estimate says that there is energy exchange, namely a beating effect, between these two modes. Notice that the sign in front of the nonlinearity does not affect the phenomena.\\
 That the principal term of $g$ be cubic is certainly not necessary but convenient for calculus. 
We will  see (see section 4) that the period of the beating depends on the nonlinearity but also on the value of all the initial actions.

The paper is organized as follows: in section 2 we prove a normal form result for the equation \eqref{cauchy} that allows to reduce the initial problem to the study of a finite dimensional classical system. The main theorem is proved in section 3. Section 4 is devoted to generalizations and comments.


\section{The normal form}
Let us expand $\psi$ and $\bar \psi$ in Fourier modes:
$$\psi(x)=\sum_{j\in \Z}\xi_j e^{ijx},\quad \bar \psi(x)=\sum_{j\in\Z} \eta_j e^{-ijx}.$$
In this Fourier setting the equation \eqref{cauchy} reads as an infinite Hamiltonian system
\be\label{hamsys}
\left\{\begin{array}{rll}
i\dot \xi_j&=j^2\xi_j+\frac{\partial P}{\partial \eta_j} &\quad j\in \Z,\\
-i\dot \eta_j&=j^2\eta_j+\frac{\partial P}{\partial \xi_j} &\quad j\in \Z  \end{array}\right.
\ee
where the perturbation term is given by
\begin{align}\begin{split}\label{P}
P(\xi,\eta)&=\pm 2\int_{S^1}\cos 2x |\psi(x)|^4 dx \\
&=\pm\sum_{\substack{j,\ell\in \Z^2\\ \mathcal M(j,\ell)=\pm 2}}\xi_{j_1}\xi_{j_2}\eta_{\ell_1}\eta_{\ell_2}.
\end{split}\end{align}
and $\Mc(j,\ell)= j_1+j_2-\ell_1-\ell_2$ denotes the momentum of the multi-index $(j,l)\in \Z^4$. In the sequel we will develop  the calculus with $\pm=+$ but all remain true, mutadis mutandi, with the minus sign. \\
 Since the regularity is not an issue in this work, we will work in the phase space ($\rho\geq 0$)
$$\Ac_\rho=\{(\xi,\eta)\in \ell^1(\Z)\times  \ell^1(\Z)\mid ||(\xi,\eta)||_\rho:=\sum_{j\in \Z}e^{\rho |j|}(|\xi_j|+|\eta_j|)<\infty\}
$$
that we endow with the canonical symplectic structure $-i\sum_j d\xi_j\wedge \eta_j$. Notice that this Fourier space corresponds to  functions $\psi(z)$ analytic on a strip $|\Im z| <\rho$  around the real axis.\\
According to this symplectic structure, the Poisson bracket between two functions f an g of $(\xi,\eta)$ is defined by
$$\lbrace{f},{g}\rbrace={-i}\sum_{j\in\mathbb{Z}}\frac{\partial{f}}{\partial{\xi_j}}\frac{\partial{g}}{\partial{\eta_j}}-
\frac{\partial{f}}{\partial{\eta_j}}\frac{\partial{g}}{\partial{\xi_j}}.$$
In particular, if $(\xi(t),\eta(t))$ is a solution of \eqref{hamsys} and $F$ is some regular Hamiltonian function, we have
$$
\frac{d}{dt}F(\xi(t),\eta(t))=\{F,H\}(\xi(t),\eta(t))
$$
where
$H=H_0+P= \sum_{j\in \Z}j^2\xi_j\eta_j +P $ is the total Hamiltonian of the system.\\
We denote by $B_\rho(r)$ the ball of radius $r$ centered at the origin in $\Ac_\rho$. In the next proposition we put the Hamiltonian $H$ in normal form up to order 4:
\begin{proposition}\label{NF}
There exists a canonical change of variable $\tau$ from $B_\rho(\e)$ into  $B_\rho(2\e)$ with $\e$ small enough such that
\begin{equation}
H\circ \tau = H_{0} +Z_{4,1}+Z_{4,2}+Z_{4,3}+R_5
\end{equation}
where
\begin{itemize}
\item[(i)]$H_{0}(\xi,\eta) = \sum_{j\in \Z}j^2\xi_j\eta_j$,
\item[(ii)] $Z_4=Z_{4,1}+Z_{4,2}+Z_{4,3}$ is the (resonant) normal form at order 4, i.e. $Z_4$ is a polynomial of order 4 satisfying $\lbrace{H}_{0},Z_4\rbrace=0$.
\item[(iii)] $Z_{4,1}(\xi,\eta)= 2(\xi_1\eta_{-1}+\xi_{-1}\eta_1)\left(2\sum_{p\in\mathbb{Z}, p\neq 1,-1}\xi_p\eta_p \ +(\xi_1\eta_1+\xi_{-1}\eta_{-1})\right)$ is the effective hamiltonian at order 4.
\item[(vi)] $Z_{4,2}(\xi,\eta)= 4\left(\xi_2\xi_{-1}\eta_{-2}\eta_1 + \xi_{-2}\xi_1\eta_2\eta_{-1}\right)$. \item[(v)] $Z_{4,3}$ contains all the term in $Z_4$ involving at most one  mode of index $1$ or $-1$.
\item[(vi)] $R_5$ is the remainder of order 5, i.e. a hamiltonian satisfying\\
 $||X_{R_5}(z)||_\rho \leq C ||z||^4_\rho$   for $z=(\xi,\eta)\in B_\rho(\e)$.
 \item[(vii)] $\tau$ is close to the identity: there exist a constant $C_\rho$ such that $||\tau(z)-z||_\rho\le{C_\rho}||z||^2_\rho$  for all $z\in B_\rho(\e).$
\end{itemize}
\end{proposition}
\proof
The proof uses the classical Birkhoff normal form procedure (see for instance \cite{Mos68} in a finite dimensional setting or \cite{Gre07} in the infinite dimensional ones). Since the free frequencies are the square of the integers, we are in a totally resonant case and we are not facing the  small denominator problems: a linear combination of integer with integer coefficients is exactly $0$ or its modulus  equals at least $1$. For convenience of the reader, we briefly recall the procedure.
Let us search $\tau$ as time one flow of $\chi$ a polynomial Hamiltonian of order 4, 
\be \label{chi}\chi = \sum_{\substack{j,\ell\in \Z^2\\ \mathcal M(j,\ell)=\pm 2}} a_{j,\ell}\   \xi_{j_1}\xi_{j_2}\eta_{\ell_1}\eta_{\ell_2}.\ee
Let $F$ be a Hamiltonian, one has by using the Taylor expansion of $F\circ \Phi^t_\chi$ between $t=0$ and $t=1$:
$$F\circ \tau = F+ \{F ,\chi \}+\frac 1 2 \int_0^1(1-t)\{\{F,\chi\},\chi\}\circ \Phi^t_\chi dt.$$
Applying this formula to $H=H_0+P$ we get
$$H\circ \tau = H_0+P+ \{H_0 ,\chi \}+\{P,\chi\}+ \frac 1 2 \int_0^1(1-t)\{\{H,\chi\},\chi\}\circ \Phi^t_\chi dt.$$
Therefore in order to obtain $H\circ \tau = H_{0} +Z_4+R_5$ we need to solve the homological equation
\be \label{homo}
\{\chi,H_0\}+Z_4=P
\ee
and then we define
\be\label{R}
R_5=\{P,\chi\}+ \frac 1 2 \int_0^1(1-t)\{\{H,\chi\},\chi\}\circ \Phi^t_\chi dt.
\ee
For $j,\ell\in \Z^2$ we define the associated divisor by
$$\Om(j,\ell)=j_1^2+j_2^2-\ell_1^2-\ell_2^2.
$$
The homological equation \ref{homo} is solved by defining
\be \label{chi}\chi = \sum_{\substack{j,\ell\in \Z^2\\ \mathcal M(j,\ell)=\pm 2, \Om(j,\ell)\neq 0}} \frac 1 {i\Om(j,\ell)} \xi_{j_1}\xi_{j_2}\eta_{\ell_1}\eta_{\ell_2}\ee
and
\be\label{Z}
Z_4 = \sum_{\substack{j,\ell\in \Z^2\\ \mathcal M(j,\ell)=\pm 2, \Om(j,\ell)= 0}} \xi_{j_1}\xi_{j_2}\eta_{\ell_1}\eta_{\ell_2}\ee
At this stage we remark that any formal polynomial $$Q=\sum_{\substack{j,l\in \Z^2\\ \mathcal M(j,\ell)=\pm 2}} a_{j,\ell}  \xi_{j_1}\xi_{j_2}\eta_{\ell_1}\eta_{\ell_2}$$ is well defined and continuous (and thus analytic) on $\Ac_\rho$ as soon as  the $a_{j\ell}$ form a bounded family. Namely, if $|a_{j,\ell}|\leq M$ for all $j,\ell\in \Z^2$ then 
$$
|Q(\xi,\eta)|\leq M||(\xi,\eta)||_0^4\leq  M||(\xi,\eta)||_\rho^4.$$
Furthermore the associated vector field is bounded (and thus smooth) from $\Ac_\rho$ to $\Ac_\rho$, namely
\begin{align}\begin{split}\label{campo}
 ||X_Q(\xi,\eta)||_\rho&=\sum_{k\in \Z} e^{\rho|k|}\left( \left|\frac{\partial Q}{\partial \xi_k}\right| +\left|\frac{\partial Q}{\partial \eta_k}\right|\right)\\
   &\leq 2M\sum_{k\in \Z} e^{\rho|k|}\sum_{\substack{j_1,\ell_1,\ell_2\in \Z\\ \mathcal M(j_1,k,\ell_1,\ell_2)=\pm 2}} \left|\xi_{j_1}\eta_{\ell_1}\eta_{\ell_2}\right|+\left|\xi_{\ell_1}\xi_{\ell_2}\eta_{j_1}\right|\\
    &\leq  2Me^{2\rho}  \sum_{j_1,\ell_1,\ell_2\in \Z} \left|\xi_{j_1}e^{\rho|j_1|}\eta_{\ell_1}e^{\rho|\ell_1|}\eta_{\ell_2}e^{\rho|\ell_2|}\right|+\left|\xi_{\ell_1}e^{\rho|\ell_1|}\xi_{\ell_2}e^{\rho|\ell_2|}\eta_{j_1}e^{\rho|j_1|}\right|\\
    &\leq 4Me^{2\rho}||(\xi,\eta)||_\rho^3.\end{split}\end{align}
    where we used, $\mathcal M(j_1,k,\ell_1,\ell_2)=\pm 2 \Rightarrow |k|\leq 2+|j_1|+|\ell_1|+|\ell_2|$.
  
Since there are no small divisors in this resonant case, $Z_4$ and $\chi$ are well defined on $\Ac_\rho$ and by construction $Z_4$ satisfies (ii). On the other hand, since $\chi$ is homogeneous of order 4,  for $\e$ sufficiently small, the time one flow generated by $\chi$  maps the ball $B_\rho(\e)$ into the  ball $B_\rho(2\e)$ and is close to the identity in the sense of assertion (vii). \\
Concerning $R_5$, by construction it is a Hamiltonian function which is of order at least 5. To obtain assertion (vi) it remains to prove that the vector field $X_{R_5}$ is smooth from $B_\rho(\e)$ into $\Ac_\rho$ in such a way we can Taylor expand $X_{R_5}$ at the origin. This is clear for the first term of \eqref{R}, $\{P,\chi\}$. For the second, notice that $\{H,\chi\}=Z_4 -P +\{P,\chi\}$ which is  a polynomial on $\Ac_\rho$ having bounded coefficients\footnote{Notice that this is not true for $H_0$.} and the same is true for $Q=\{\{H,\chi\},\chi\}$. Therefore, in view of the previous paragraph,  $X_Q$ is smooth. Now, since for $\e$ small enough $\Phi^t_\chi$ maps smoothly the ball $B_\rho(\e)$ into the  ball $B_\rho(2\e)$  for all $0\leq t\leq 1$, we conclude that $\int_0^1(1-t)\{\{H,\chi\},\chi\}\circ \Phi^t_\chi dt$ has a smooth vector field.

We are now interested in describing more explicitly the terms appearing in the  expression \eqref{Z}  on  base of the number of times that $\xi_{\pm{1}}$ or
 $\eta_{\pm{1}}$ appear as a factor. Let us define the resonant set
 \begin{equation}
{\mathcal R}=\lbrace(j_1,j_2,j_3,j_4)\in\mathbb{Z}^4 \mid j_1+j_2-j_3-j_4=\pm{2} \;\; {\rm and } \;\; j_1^2+j_2^2-j_3^2-j_4^2=0\rbrace
\end{equation}
in such a way \eqref{Z} reads
\begin{equation}\label{expz4}
Z_4= \sum_{(j_1,j_2,j_3,j_4)\in{\mathcal R}}\xi_{j1}\xi_{j2}\eta_{j3}\eta_{j4}.
\end{equation}
Note that it is enough to deal with the case $j_1+j_2-j_3-j_4=2$ because the other case   $j_1+j_2-j_3-j_4=-2$ can be
 obtained by changing the signs of  $j_1,j_2,j_3,j_4$.
 Thus let us study two cases:  $j_1=1$ or $j_1=-1$.

 In the first case we have $j_2=j_3+j_4+1$ and therefore the quadratic condition implies $(1+j_3)(1+j_4)=0$.  Thus we conclude that $(j_1,j_2,j_3,j_4)=(1,p,-1,p)$ or $(j_1,j_2,j_3,j_4)=(1,p,p,-1)$. Thus  when $j_1=1$, the only terms appearing in the expression for $Z_4$ in Eq. (\ref{expz4}) must have the
 form $\xi_1,\xi_p,\eta_{-1}\eta_p$ or $\xi_1,\xi_p,\eta_p,\eta_{-1}$.

In the second case, $j=-1$, we have $j_2=j_3+j_4+3$ which in turn implies $5+3(j_3+j_4)+j_3j_4=0$. By looking at the graph of the function
$f(x)=-\frac{5+3x}{3+x}$ we find that the only acceptable cases appearing in Eq. (\ref{expz4}) are $(j_1,j_2,j_3,j_4)=(-1,1,-1,-1),\;
(-1,2,-2,1),\; (-1,2,1,-2),\; (-1,-8,-4,-7),\\ (-1,-8,-7,-4)\;$, or $(-1,-7,-5,-5)$.

Noting that the resonant set ${\mathcal R}$ is invariant under permutations of $j_1$ with $j_2$ or $j_3$ with $j_4$ we obtain that
\begin{eqnarray*}
  Z_4  &=& \left(4\sum_{p\in\mathbb{Z}}\xi_p\eta_p-2(\xi_1\eta_1+\xi_{-1}\eta_{-1})\right)\left(\xi_1\eta_{-1}+\xi_{-1}\eta_1\right) \\  
  && + 4\left(\xi_2\xi_{-1}\eta_{-2}\eta_1 + \xi_{-2}\xi_1\eta_2\eta_{-1}\right) \nonumber\\
   && + 4\left(\xi_{-1}\xi_{-8}\eta_{-7}\eta_{-4} + \xi_1\xi_8\eta_7\eta_4\right) + 2\left(\xi_{-1}\xi_{-7}\eta_{-5}\eta_{-5}+\xi_{1}\xi_{7}\eta_{5}\eta_{5}\right) 
   + \tilde{Z_4}\end{eqnarray*}
where $\tilde{Z_4}$ is equal to the sum of terms of the form  $\xi_{j1}\xi_{j2}\eta_{j3}\eta_{j4}$  with $(j_1,j_2,j_3,j_4)\in{\mathcal R}$ and satisfying 
the condition $|j_k|\neq{1}$, for all $k=1,2,3,4$. 

\endproof

\section{Dynamical consequences}\label{dyn}
 We denote by $I_p=\xi_p\eta_p$, $p\in \Z$, the actions, by $J_p=I_p+I_{-p}$ for $p\in \N\setminus\{0\}$ and $J_0=I_0$, the generalized actions  and by $J=\sum_{p\in \Z}I_p$. Notice that, when the initial condition $(\xi^0,\eta^0)$ of the Hamiltonian system \eqref{hamsys} satisfies $\eta^0=\bar \xi^0$ (and this is actually the case when $\xi^0$ and $\eta^0$ are the sequences of Fourier coefficients of respectively $\psi_0$, $\bar \psi_0$), this reality property is conserved i.e. $\eta(t)=\bar \xi(t)$ for all $t$. As a consequence all the quantities $I_p$, $J_p$ and $J$ are real and positive.
 
   To begin with we establish that apart from $I_1$ and $I_{-1}$ the others actions are almost constant:
\begin{lemma}\label{lem}
Let $\psi(t,\cdot)=\sum_{k\in\Z} \xi_k(t) e^{ikx}$ be the solution of \eqref{cauchy} then for $\e$ small enough and $|t|\leq \e^{-\frac 9 4}$,
$$|J(t)-J(0)|\leq \e^{5/2} \mbox{ and } J_p(t)\leq \e^3 \mbox{ for } p\in \N\setminus\{1\} .$$
As a consequence
$$
 I_1(t),\  I_{-1}(t)\leq 4\e^2
\mbox{ while  }  I_p(t)\leq \e^3 \mbox{ for } p\in \N\setminus\{1\} .
$$
\end{lemma}
\proof
Using Proposition \ref{NF} we have for all $p\in \N$
 $$\dot J_p=\{J_p,H\}=\{J_p,Z_{4,3}\}+\{J_p,R\}$$
 since by  elementary calculations $\{J_p,Z_{4,1}\}=\{J_p,Z_{4,2}\}=0$. In particular
 $$\dot J=\{J,H\}=\{J,Z_{4,3}\}+\{J,R\}.$$
 Let $T_\e$ be the maximal time such that for all $|t|\leq T_\e$
 $$|J(t)-J(0)|\leq \e^{5/2} \mbox{ and }  J_p(t)\leq \e^3,\mbox{ for } p\in \N\setminus\{1\}.$$
 Since $J(0)=\e^2$, we get for $\e$ small enough and $|t|\leq T_\e$
 $$
 |\xi_1(t)|,\  |\xi_{-1}(t)|,\ |\eta_1(t)|,\  |\eta_{-1}(t)|\leq 2\e
\mbox{ while  }  |\xi_p(t)|,\ |\eta_p(t)|\leq \e^{3/2} \mbox{ for } p\neq 1 .
$$
Therefore, from the definition of $Z_{4,3}$ and $R$ we deduce that for $|t|\leq T_\e$
\begin{align*}\{J_p,R\}&=O(\e^{\frac{11}2}),\  \{J_p,Z_{4,3}\}=O(\e^{\frac{11}2})\mbox{ for } p\neq 1,\\
 \{J,R\}&=O(\e^{5}),\ \{J,Z_{4,3}\}=O(\e^{\frac{11}2}),\end{align*}
and thus there exists $C>0$ such that for any $|t|\leq T_\e$
$$\ |\dot J(t)|\leq C\e^{5} \mbox{ and } |\dot J_p|\leq C\e^{\frac {11} 2} \mbox{ for } p\neq 1.$$
Taking into account that $J_p(0)=0$ for $p\neq 1$, we obtain 
$$|J(t)-J(0)|\leq C|t|\e^5,  \mbox{ and } | J_p(t)|\leq C|t|\e^{\frac {11} 2} \mbox{ for } p\neq 1
$$ for all $|t|\leq T_\e$ and we conclude by a classical bootstrap argument that $T_\e\geq C^{-1}\e^{-5/2} \geq \e^{-9/4}$ for $\e$ small enough.
\endproof
In order to prove theorem \ref{main} let us define some quadratic Hamiltonian functions ($p\in \Z)$:
\begin{align*}\label{quadra}M_p&:=\xi_p\eta_{p}-\xi_{-p}\eta_{-p}, \quad J_p:=\xi_p\eta_{p}+\xi_{-p}\eta_{-p},\\ L_p&:=\imath(\xi_p\eta_{-p}-\xi_{-p}\eta_{p}), \quad K_p:=\xi_p\eta_{-p}+\xi_{-p}\eta_{p}.
\end{align*}
One computes
\begin{align*}\dot M_1&=\{M_1,H\}=\{M_1,Z_{4,1}+Z_{4,2}\}\ +\ \{M_1,Z_{4,3}+R\}\\
&=2JL_1+L_1K_2-K_1L_2\ +\ \{M_1,Z_{4,3}+R\}
\end{align*}
and using Lemma \ref{lem} we get that for $|t|\leq  \e^{-9/4}$,
$$\dot M_1= 2J(0)L_1+O(\e^{9/2})$$
and in the same way we verify
$$\dot L_1= -2J(0)M_1+O(\e^{9/2}).$$
We can now compute the solution of the associated linear ODE 
\begin{equation*} 
\left\{\begin{array}{ll}\dot M_1&= 2J(0)L_1\\
\dot L_1&= -2J(0)M_1
\end{array}\right.
\end{equation*}
to conclude that for $ t\leq C\e^{-9/4}$
$$
M_1(t) =M_1(0)\cos {2J(0)} t+L_1(0) \sin {2J(0)} t + O(\e^{9/4}).
$$
In Theorem \ref{main} we have chosen $\psi_0=\e(\cos x+\sin x)$ which corresponds to $\xi_1=\eta_{-1}=\frac{1-i}{2}\e$ and $\eta_1=\xi_{-1}=\frac{1+i}{2}\e$. Therefore $J(0)=J_1(0)=\e^2$, $M_1(0)=0$ and $L_1(0)=\e^2$ which leads to the desired result. \\
We finally remark that choosing the minus sign in front of the nonlinearity will lead to the following linear system
\begin{equation*} 
\left\{\begin{array}{ll}\dot M_1&= -2J(0)L_1\\
\dot L_1&= 2J(0)M_1
\end{array}\right.
\end{equation*}
which again gives the desired result.
\section{Generalizations and comments}
\bi
\item The same result remains true when we add a higher order term to the nonlinearity, i.e. considering the equation
$$i\psi_t= -\psi_{xx}\pm 2\cos 2x \ |\psi|^2\psi +O(|\psi|^4),\quad x\in S^1,\ t\in \R.$$
\item We can prove a similar result when changing the nonlinearity in such a way we still privilege the modes 1 and -1. The game is to conserve an effective Hamiltonian at order 4 (see Proposition \ref{NF}). For instance $2\cos 2x $ can be replace by $a\cos 2x +b\sin 2x$ but not by $\cos 4x$ which generates an effective Hamiltonian only at order 6. We can also choose to privilege another couple of modes $p$ and $-p$ choosing $a\cos 2px +b\sin 2px$. In that case we have also to adapt the initial datum. 
\item If we choose a non linearity that does not depend on $x$, for instance the standard cubic nonlinearity $g=\pm |\psi|^2\psi$, then we can prove that there is no beating effect between any modes for $|t|\leq \e^{-3}$ since $Z_4$ only depends on the actions in that case (independently of the sign in front of the nonlinearity). 
\item We can also change the initial datum. Remark that if you chose $\psi_0=\e \cos x$ or $\psi_0=\e \sin x$, no beating effect  appears at order 4 since in both cases $M_1(0)=L_1(0)=0$.  
\item The beating frequency, namely $2J(0)$, depends on all the modes initially excited. For instance if $\psi_0=\e(\cos x+\sin x) + \e^2\cos qx$ ($q\neq 1$)  then we can still prove  that there is no energy  exchanges between the mode $p$, for $|p|\neq 1$, and modes $1$ and $-1$, that there is the same beating effect between modes $1$ and $-1$, nevertheless the beating frequency is slightly changed: $2J(0)=2\e^2+\e^4$. 
\item As stated in the introduction, when adding a linear potential –a multiplicative one or a convolution one–  we can choose the potential in order to avoid resonances between the different blocks of modes $p$ and  $-p$ (see \cite{BG06} for multiplicative potentials or \cite{Gre07} for convolution potentials). In that case the same result can be proved and actually in an easier way since we avoid exchanges between the blocks for arbitrary long time.

\ei

\bibliographystyle{amsalpha}

\begin{thebibliography}{99}

  
\bibitem[BG06]{BG06}
D.~Bambusi and B.~Gr{\'e}bert, \emph{Birkhoff normal form for {PDE}s with tame
  modulus}, Duke Math. J. \textbf{135} (2006), 507--567.
  
\bibitem[Bou99]{Bou99}
  J. Bourgain, ,
   \emph{Global solutions of nonlinear Schr\"odinger equations},
   {American Mathematical Society Colloquium Publications},
  { \bf 46},
  Amer. Math. Soc. Providence RI, {1999}.


  
\bibitem[Gr{\'e}07]{Gre07}
B.Gr{\'e}bert, \emph{Birkhoff normal form and {H}amiltonian
  {PDE}s}, Partial differential equations and applications, S\'emin. Congr.,
  vol.~15, Soc. Math. France, Paris, 2007, pp.~1--46. 

\bibitem[Mos68]{Mos68}
J.~Moser, \emph{Lectures on hamiltonian systems}, Mem. Amer. Math. Soc.
  \textbf{81} (1968), 1--60.

\end{thebibliography}

{\bf Beno\^{i}t Gr\'ebert}

{\it Laboratoire de Math\'ematique Jean Leray UMR 6629,

Universit\'e de Nantes,
2, rue de la Houssini\`ere,

44322 Nantes Cedex 3, France}

E-mail: \quad
\parbox{5cm} {\tt benoit.grebert@univ-nantes.fr}

\bigskip

{\bf Carlos Villegas-Blas}

{\it Universidad Nacional Autonoma de México,

Instituto de Matem\'aticas,
Unidad Cuernavaca}

E-mail: \quad
\parbox{5cm} {\tt villegas@matcuer.unam.mx}

\end{document}